\documentclass[11pt,a4paper]{article}
\usepackage[utf8]{inputenc}
\usepackage[english]{babel}
\usepackage{amssymb, amsmath, enumerate}
\usepackage{graphicx}
\usepackage{epstopdf}
\usepackage{subcaption}
\usepackage{authblk}
\usepackage{textcomp}

\newcommand{\qed}{\nobreak \ifvmode \relax \else
	\ifdim\lastskip<1.5em \hskip-\lastskip
	\hskip1.5em plus0em minus0.5em \fi \nobreak
	\vrule height0.75em width0.5em depth0.25em\fi}
\begin{document}
	\title{The problem of optimal location of production points and distribution points in the vertices of the  transportation network as an investment project}
	\author[1]{Malafeyev O.A.\thanks{o.malafeev@spbu.ru}}
	\author[1]{Onishenko V.E.\thanks{st034043@student.spbu.ru}}
	\affil[1]{Saint-Petersburg State University,  Russia}
	\date{}
	\maketitle
	\begin{abstract}
	The article deals with the problem of optimal location of production and distribution points in the vertices of the transportation network considered as an investment project. On the network set production of raw materials points and demand points.
	
	\end{abstract}
	
	

	\section{Introduction}

Facility location problems study how to best locate facilities under
the assumption that consumers will go to the facility that is most
profitable to them. Hotelling (1929) introduced an influential competitive facility location model where each of two players select a
location in a linear segment and price of consumption is constant. Economic agents (players, owners), in this model, is to maximize their demand. It's model a non-cooperative game with equilibrium.
 Facility location games have several applications related rumor dissemination, seeding in social networks and product differentiation models (Bharathi et al. 2007, Kostka et al. 2008). The following articles consider spatial competition in a continuous space with a Nash equilibrium. ( Eaton and Lipsey 1976, Graitson 1982, D’Aspremont et al. 1983, de Palma et al. 1985, Gabszewicz and Thisse 1986). In the following articles, a model of the circular city of Salopa is considered (Salop 1979, Eiselt and Laporte 1989, Eiselt et al. 1993, Eiselt and Laporte 1996, Plastria 2001, ReVelle and Eiselt 2005, Smith et al. 2009, Dasci 2011, Kress and Pesch 2012).
The following articles consider the possibility of the existence of Nash equilibrium for complex network structures( D¨urr and Thang 2007). Nash equilibrium, Stackelberg games and models where agents locate facilities sequentially are considered in several articles
( Hay 1976, Prescott and Visscher 1977, Drezner 1982, ReVelle 1986, Hakimi 1983, Hakimi 1986, Leonardi and Tadei 1984, Drezner and Drezner 1996, Drezner 1998, Leonardi and Tadei 1984, Drezner and Drezner 1996) applying a gravity rule and defining a sphere of influence ( Drezner 1995a, Drezner et al. 2002, Drezner and Drezner 2004, ReVelle 1986, Drezner et al. 2011 ).
Co-operative games and competition with prices (Vetta 2002, Mirrokni and Vetta 2004, Jain and Mahdian 2007).
	
The following problem is considered. In some vertices of the network there are goods demand points and production points of raw materials. On the set of network edges, the function of transport-corruption costs is defined, which denotes the cost of moving along the edges of the network between the vertices of the network. All vertices of the network are connected, that is, from any vertex of the network there
	exists a possibly non-unique path 
	to any other vertex of the network.
	
	For each demand point for a product, the quantity of goods that it wishes to buy is given. The amount of expenses for meeting the demand is equal to: the amount of costs for the purchase of goods, the amount of costs for moving between the vertices of the network to the goods distribution points. It is assumed that each demand point for goods chooses the distribution point of goods to meet their demand for goods, based on minimizing their amount of costs to meet their demand for the goods.

Owners of goods distribution points and goods production points wish to location their goods production points and goods distribution points at the vertices of the network in the most profitable way for them in terms of maximizing the revenue from selling this product. 
Thus arises the optimal location problem
of goods distribution points and goods production points on a given network, in accordance with selected principle of optimality.

Papers related to the theme of this article are [1 - 35]. 
	
\section{Formalization of the problem of optimal location of  goods production points and goods distribution points at the vertices of the network}

Let a network with a capacity $ (N, k) $ be given on a plane, where $ N $ is a finite set of nodes, $ k $ is a capacity function that maps to each edge of the network $ (x, y), (x, y) \in N $ is a nonnegative number $ k (x, y) \geq0 $. Let also in some vertices of the network making up the set $ L = (l_1, ..., l_h) $, there are production points of raw materials necessary for the production of goods at the production points  of the goods, and at the vertices composing the set $ K = (k_1 ,. .., k_s) $ are the demand points for the goods.We denote the set of goods production points by $ M = (m_1 ... m_n) $. The set of goods distribution points will be denoted by $ W = (w_1 ... w_n) $. On the set $ E $ of edges of the network $ (N, k) $, the transport-corruption cost function $ C: E \rightarrow R_1, C (x, y) \geq0 $ is given. It is assumed that the vertices of the network $ (N, k) $ are connected to each other, that is, from any vertex of the network there 
exists a path to any other vertex of the network, or perhaps multiple.

Owners of goods production points and goods distribution points $ R = (R_1, ..., R_n) $ assign prices for the goods they produce $ P = (p_1 ... p_n) $. In doing so, they tend to maximize their profits. We believe that the demand points for the goods satisfy their demand at the point of distribution of the goods, where the total costs are minimal. Net income of the owners is the sum of the funds received by them from the point of demand for the goods, minus the cost of production of goods in locations of production of goods, the costs of distribution of goods to goods distribution points, the costs for the purchase of raw materials in the points production of raw materials and transport and corruption costs.

$ TC (R_i) $ - is the net income of the owner $ R_i $.

$$
TC(R_i) = P(w_i)D(w_i) - ( C(R_i)+PW(R_i)+PM(R_i)+PL(R_i) )
$$
$ P (w_i) $ - is the price assigned at the goods distribution points $ w_i $ for the goods.\\
$ D (w_i) $ - is the demand at the distribution point of the commodity $ w_i $, which consists of the demand of all points of demand for the goods that satisfy their demand at the distribution point of goods $ w_i $.\\
$ C (R_i) $ - all transport-corruption costs of the owner of goods production points and goods distribution points $ R_i $.\\
$ PW (R_i) $ - costs of the owner $ R_i $ of goods distribution points and points of production of goods for the distribution of goods.\\
$ PM (R_i) $ - costs of the owner $ R_i $ of goods distribution points and points of production of goods for the production of goods.\\
$ PL (R_i) $ - the total cost of raw materials for the owner $ R_i $.\\

The net income of each owner of the goods distribution points and goods production points is represented by a table, by columns of which we indicate the net incomes of each owner $ TC (R_i) $, and by rows the possible locations of goods production points and goods distribution points at the vertices of the network and the price $ p_i $ for unit of output.

\begin{center}
    \fontsize{9pt}{9pt}\selectfont

\begin{tabular}{l||lll|}

&$TC(R_1)$&{...} \\
\hline
\hline

$(p_1,m_1,w_1)...(p_1,m_{n-1},w_{n-1}),(p_1,m_n,w_t)$  & $ TC(R_1(p_1,m_1,w_1)) $& {...}  \\
$(p_1,m_1,w_1)...(p_1,m_{n},w_{t}),(p_1,m_{n-1},w_{t-1}) $&$ TC(R_1(p_1,m_1,w_1))$ & {...} \\
{...} & {...} & {...} \\ 
$(p_2,m_1,w_1)...(p_1,m_{n-1},w_{t-1}),(p_1,m_n,w_t)$  &$ TC(R_1(p_2,m_1,w_1)) $& {...} \\
{...} & {...}& {...}  \\ 
$(p_2,m_1,w_1)...(p_2,m_{n-1},w_{t-1}),(p_2,m_n,w_t)$ & $TC(R_1(p_2,m_1,w_1)) $& {...} \\
{...} & {...} & {...} \\ 
$(p_l,m_n,w_n)...(p_l,m_{2},w_{2}),(p_l,m_1,w_1) $&$ TC(R_1(p_l,m_n,w_t))$ & {...}  \\

\end{tabular}
\end{center}

\qquad \qquad

Owners of goods distribution points and
goods production points agree among themselves on the location of their own
points of production of goods and distribution points at the vertices of the network on the basis of some optimality principle. For example, a compromise solution or a Cournot-Nash equilibrium, algorithms for finding
which are described below.

The compromise set is defined by the following formula:
$$
C_H= \{ x\in X |\max_i(M_i
-H_i(x)) \le \max_i (M_i -H_i(x')) \forall  x'\in X \}.
$$

Otherwise, the compromise set can be defined as follows

Let $ X $ be the set of admissible solutions

$ H_i: X \rightarrow R_1 $
  - the income function of the agent $ i $ on the solution set.

Let $ M $ = $ (M_1 ... M_n) $

$$
M_i = \max_ {x \in X} H_i (x)
$$

$ M $ is an ideal vector

We choose a fixed $ x \in X $ and calculate the residual vector for it.

$$
\Delta (x) = \{M_i - H_i (x) \} _ 1 ^ n
$$

For each $ x \in X $

$$
\max_ {i \in I} \{M_i-H_i (x) \} _ 1 ^ n
$$

The compromise set $ C_X ^ n $ is
$$
arg
\min_ {x \in X} \max_ {i} \delta (x) = C_X ^ n
$$

\qquad

Algorithm for finding a compromise solution.

\qquad

\emph{Step 1.} Let $ Z $ be the set of admissible solutions $ z_s \in Z $,

where $ z_s $ = $ ((p_1 ^ {i_1} m_1 ^ {j_1} w_1 ^ {q_1}), ..., (p_l ^ {i_l} m_n ^ {j_n} w_t ^ {q_t})) $,

$i_1 = \overline {1, \overline {i_1}}$,
$j_1 = \overline {1, \overline {j_1}}$, $w_1 = \overline {1, \overline {w_1}}$, $i_l = \overline {1, \overline {i_l}}$, $j_n= \overline {1, \overline {j_n}}$,
$w_q = \overline {1, \overline {w_q}}$,

- these are the possible vertices of the location of goods production points and goods distribution points at specified prices, and $ s $ is the index by which all permissible solutions are renumbered and $ s = 1, ..., \overline {s}$.\\
TC ($ R_1 $ ($ z_s $)), ..., TC ($ R_n $ ($ z_s $)) is the net income of the i owner of the goods production points and distribution points for a given $ z_s $. Let us construct an ideal vector
$ M $ = [$ M_1 $, ..., $ M_n $] = ($ max $ ($ TC $ ($ R_1 $ ($ z_s $))), ..., $ max $ ($ TC $ ($ R_n $ ($ z_s $)))) is the maximum income of the j-th owner of the goods production points and distribution points.

\emph {Step 2.} For each feasible solution $ z_s, s = 1, ..., \overline {s} $, we form the residual matrix, which in this case will have the form

$ A * = (M_1-TC (R_1 (z_s))), ..., (M_n-TC (R_n (z_s))) = (\alpha_s ... \beta_s) $,

where $ s = 1 ... \overline {s} $.

\emph {Step 3.} In the resulting residual matrix, for each solution $ z_s, s = 1, ..., \overline {s} $ from the columns $ (\alpha_s ... \beta_s) $, we select the maximum value $ \gamma_s = \max\limits_s (\alpha_s ... \beta_s) $.

\emph {Step 4.} We choose the minimal of these maximal solutions $ min = \gamma_s $, which will be a compromise solution.

\qquad

Algorithm for finding the Cournot-Nash equilibrium.

\qquad

The Cournot-Nash equilibrium is the type of solution of a game of two or more players in which neither of the participants can increase his winnings by changing his decision if other participants do not change their decisions. Let's make a table, according to the columns and rows of which
we can get
the possible location of the goods production points and distribution points at the vertices of the network at the indicated prices. As a payoff function, the net income of each owner of the goods production points and goods distribution points is considered.

\emph {Step 1.} Find the maximum value of winning by rows.

\emph {Step 2.} Find the maximum value of winning by columns.

\emph {Step 3.} We consider the intersection of the obtained maxima. Where 

the intersection exists, this is the Cournot-Nash equilibrium.

\qquad

Algorithm for solving the problem of optimal location of goods production points and goods distribution points at the vertices of the network.

\qquad

1. Set a price $ P = (p_1 ... p_n) $ per unit of the produced goods for each vertex of the possible location of goods distribution points $ W = (w_1 ... w_n) $;

2. Calculate the total cost goods for goods demand points, taking into account the costs of purchasing the goods and the transport-corruption costs of moving to the goods distribution points. Using the Floyd algorithm, to find the shortest path between the vertices of a graph,
to count transport-corruption costs. For each of demand point for goods, determine the optimal distribution point of goods to meet demand;

3. Calculate the amount of demand for goods in the $ w_i $ distribution point of goods, which equals the sum of the demand for the goods of all goods demand points that satisfy their demand in the $ w_i $ distribution point of goods;

4. Calculate all the transport-corruption costs of $ C (m_n, w_n) $ for each possible location of goods production points $ M = (m_1 ... m_n) $ and goods distribution points $ W = (w_1 ... w_n) $, using Floyd's algorithm for finding the shortest path between the vertices of the graph;

5. Calculate the net income of $ TC (R_i (p_l,m_n, w_t)) $ of each owner of goods distribution points and goods production points;

6. Using the algorithm for finding a compromise solution or the algorithm for finding the Cournot-Nash equilibrium, we find the optimal location of the goods production points and goods distribution points;

\section{Example}

Let $ (N, K) $ be a network with a capacity, where $ N $ is the vertex set and the capacity function k that associates with each edge $ (x, y) $ a nonnegative number $ k (x, y) $ containing 8 vertices and 7 edges. For each network edge transport-corruption costs are set.

The function of transport-corruption costs of the $ (N, k) $ network is represented by the matrix $ C $, where columns and rows indicate transport-corruption costs for each edge of the network:

\qquad 

\qquad \qquad \qquad 
\begin{tabular}{c||rrrrrrrr|}
&1&2&3&4&5&6&7&8 \\
\hline
\hline

1&0&1&$\infty$&$\infty$&$\infty$&$\infty$&$\infty$&$\infty$ \\
2&1&0&2&2&2&$\infty$&$\infty$&$\infty$ \\
3&$\infty$&2&0&$\infty$&$\infty$&$\infty$&$\infty$&$\infty$ \\
4&$\infty$&2&$\infty$&0&$\infty$&$\infty$&$\infty$&$\infty$ \\
5&$\infty$&2&$\infty$&$\infty$&0&3&$\infty$&$\infty$ \\
6&$\infty$&$\infty$&$\infty$&$\infty$&3&0&3&3 \\
7&$\infty$&$\infty$&$\infty$&$\infty$&$\infty$&3&0&$\infty$ \\
8&$\infty$&$\infty$&$\infty$&$\infty$&$\infty$&3&$\infty$&0 \\

\hline
\end{tabular}

\qquad 

$ L = (l_1 ^ {1}) $ - location of the production point of raw materials

$ P (l_1 ^ {1}) = 1 $ - the cost of a raw material unit

$ K = (k_1 ^ 4, k_2 ^ 5, k_3 ^ 8) $ - location goods demand points

$ D = (10,10,10) $ - the value of demand in the demand points

$ W = (w_1 ^ 2, w_2 ^ 6) $ - possible location of goods distribution points

$ PW = (10,10) $ - value of distribution costs

$ M = (m_1 ^ 3, m_2 ^ {7}) $ - possible location of goods production points

$ PM = (40,40) $ - the cost of production costs

$ P = (p_1) = (10) $ - possible prices for the goods

We use the Floyd algorithm to find the shortest path between the vertices of the network and the calculated transport-corruption costs for the goods demand points.

\qquad \qquad \qquad \qquad \qquad \qquad 
\begin{tabular}{c||rrrr|}

&$w_1$&$w_2$ \\
\hline
\hline
$k_1$&2&7 \\
$k_2$&2&3 \\
$k_3$&8&3 \\

\hline
\end{tabular}

\qquad \qquad 

The income of distribution points for given prices is $(p_1,w_1; p_1,w_2)$.

\qquad \qquad 

\begin{center}
 
\begin{tabular}{c||cccc|}

&$P(w_i)D(w_i)$ \\
\hline
\hline
$w_1$& 200 \\
$w_2$& 100 \\

\hline
\end{tabular}
\end{center}

We will calculate the costs of each combination of goods distribution points and points of production
goods:

\begin{center}
    \begin{tabular}{c||cccc|}

&$C(w_1,m_1)+PW+PM$ \\
\hline
\hline
$(w_1,m_1)$&75 \\
$(w_1,m_2)$&87 \\
$(w_2,m_1)$&70 \\
$(w_2,m_2)$&72 \\

\hline
\end{tabular}

\end{center}

We calculate the net income $ TC (w_i, m_i) $

\begin{center}
    \begin{tabular}{c||cccc|}

&$TC(w_i,m_i)$ \\
\hline
\hline
$(w_1,m_1)$&125 \\
$(w_1,m_2)$&113 \\
$(w_2,m_1)$&30 \\
$(w_2,m_2)$&28 \\

\hline
\end{tabular}

\end{center}

The Cournot-Nash equilibrium.

\begin{center}
\begin{tabular}{c||cccc|}

&  $TC(R_2(p_1^2,w_1^2,m_1^2))$  &  $TC(R_2(p_1^2,w_1^2,m_2^2))$  \\
\hline
\hline
$TC(R_1(p_1^1,w_1^1,m_1^1))$&    (125,0) &  (125,0)            \\
$TC(R_1(p_1^1,w_1^1,m_2^1))$&    (113,0 ) & (113,0 )            \\
$TC(R_1(p_1^1,w_2^1,m_1^1))$&    (30,0  ) &(30,113 )            \\
$TC(R_1(p_1^1,w_2^1,m_2^1))$&   ( 28,125 )& (28,0  )          \\

\hline
\end{tabular}

\end{center}

\begin{center}

\begin{tabular}{c||cccc|}

&    $TC(R_2(p_1^2,w_2^2,m_1^2))$  &  $TC(R_2(p_1^2,w_2^2,m_2^2))$ \\
\hline
\hline
$TC(R_1(p_1^1,w_1^1,m_1^1))$&    (125,0 )& (125,28)         \\
$TC(R_1(p_1^1,w_1^1,m_2^1))$&     (113,30) &(113,0 )         \\
$TC(R_1(p_1^1,w_2^1,m_1^1))$&     (30,0 ) & (30,0 )           \\
$TC(R_1(p_1^1,w_2^1,m_2^1))$&     (28,0 )&  (28,0 )           \\

\hline
\end{tabular}

\end{center}

\qquad

The equilibrium solution according to this table is $ (125, 28) $. The location of the  goods production point of the first owner of the production and goods distribution points $ (R_1) $ is at vertex 3, and the goods distribution point is at vertex 2, at the selected price $ p_1 $. The location of the production point of the goods of the second owner of the production and distribution points $ (R_2) $ is at the vertex 7, and the item distribution point at the vertex 6, is at the selected price $ p_1 $.

\section{Conclusion}

Using the work model of location in the vertices of the network of goods production points  and goods distribution points, given the location of points of production of raw materials and goods demand points, one can find the optimal location of goods production points and distribution points at the vertices of the network. The presented model can be used in practice. As an example, we can consider vertically integrated engineering companies,
that wish to locate their engineering plant and storage in different cities located at different distances from each other.

%\newline
	\section{Acknowledgements}
	The work is partly supported by work RFBR No. 18-01-00796.

\end{document}